\documentclass[12pt]{amsart}

\usepackage{amsmath,amsthm, amssymb}
\usepackage[T1]{fontenc}
\usepackage{latexsym}
\usepackage{subfigure, epsfig, epsf}
\usepackage{color,fancyhdr,lastpage,graphicx}
\usepackage{xspace}
\usepackage{wasysym}
\usepackage{setspace}
\usepackage{bbm}
\usepackage{stmaryrd}

\usepackage{pstricks}
\usepackage{pst-node}
\usepackage{pst-tree}
\usepackage{todonotes}

\newcommand{\RR}{\mathbb{R}}        
 \newcommand{\OO}{\mathbb{O}}

 \newcommand{\mcC}{\mathcal{C}}   
\newcommand{\mcH}{\mathcal{H}} \newcommand{\mcV}{\mathcal{V}}

  \newcommand{\be}{{\bf e}}

    \newcommand{\bfX}{{\bf X}}  \newcommand{\bfY}{\bf Y}   \newcommand{\bfB}{\bf B} 
\newcommand{\bfZ}{\bf Z}    \newcommand{\bfM}{\bf M}    \newcommand{\bfF}{\bf F}    
\newcommand{\bfN}{\bf N}   \newcommand{\bfH}{\bf H}

\newcommand{\ep}{\epsilon}
\renewcommand{\leq}{\leqslant}
\renewcommand{\geq}{\geqslant}

\newcommand{\ssk}{\smallskip}

\newtheorem{thm}{\hspace{-0.15cm}  {\sc Theorem} }
\newtheorem{prop}[thm]{\hspace{-0.15cm} {\sc Proposition}}
\newtheorem{defn}[thm]{ \hspace{-0.3cm} {\sc Definition}}

\numberwithin{equation}{section} 

\newenvironment{Dem}{%
    \begin{list}{\hspace{0.5cm}{\sc Proof --}}{%
        \setlength{\topsep}{0pt}%
        \setlength{\leftmargin}{0pt}%
        \setlength{\rightmargin}{0pt}%
        \setlength{\listparindent}{0pt}%
        \setlength{\itemindent}{0pt}%
        \setlength{\parsep}{0pt}%
        \addtolength{\leftmargin}{20pt}%
        \addtolength{\rightmargin}{0pt}%
    } \item }{\hfill{\space $\rhd$}\end{list}\smallskip}

\title{Rough integrators on Banach manifolds}
\date{\today}
\author{I. Bailleul}
\address{IRMAR, 263 Avenue du General Leclerc, 35042 RENNES, France}
\email{ismael.bailleul@univ-rennes1.fr}
\thanks{This research was partially supported by the program ANR "Retour Post-Doctorants", under the contract ANR 11-PDOC-0025. The author thanks the U.B.O. for their hospitality.}

\begin{document}

\maketitle

\begin{abstract}
We introduce a notion of $p$-rough integrator on any Banach manifolds, for any $p\geq 1$, which plays the role of weak geometric H\"older $p$-rough paths in the usual Banach space setting. The awaited results on rough differential equations driven by such objects are proved, and a canonical representation is given if the manifold is equipped with a connection.
\end{abstract}

\section{Introduction}
\label{SectionIntroduction}

One of the most basic and illustrative example of a rough path is associated with the following 2-dimensional spinning signal
$$
h^n_t = \frac{1}{n}\big(\cos(n^2t),\sin(n^2t)\big), \quad 0\leq t\leq 1.
$$
Although $h^n$ converges uniformly to $0$, this path sweeps an area $t$ in any time interval $[0,t]$, independantly of $n$, suggesting that this sequence of paths may be really different from the zero signal. As a matter of fact, given two vector fields $V_1, V_2$ on $\RR^2$, it is remarkable that the solution path $x^n_\bullet$ to the ordinary controlled differential equation
$$
\dot x^n_t = \big(\dot h^n_t\big)^iV_i\big(x^n_t\big)
$$
converges uniformly to the solution path $x_\bullet$ to the ordinary differential equation
$$
\dot x_t = \big[V_1,V_2\big](x_t).
$$
This is due to the fact that the canonical lift of $h^n_\bullet$ as a weak geometric H\"older $p$-rough path, for any $2<p<3$, converges in a rough paths sense to the pure area rough path with area process $(t-s)\begin{pmatrix}
0 & 1 \\
-1 & 0
\end{pmatrix}$.

\ssk

This elementary example makes it clear that, when understood as controls, smooth paths may have a richer structure than expected at first sight, as the above constant path with non-trivial area process shows. Smooth paths have canonical lifts as rough paths; a smooth path in $\RR^\ell$ whose rough path lift is non-canonical will be called a {\it spinning path}. While this may not be obvious, we shall see in this note that it actually makes perfect sense to talk of a spinning manifold-valued path, as an example of a manifold valued 'rough path'. One should keep in mind, however, that the classical notion of rough path only becomes really interesting when understood as a control in some differential equation. Does it make sense to think of a manifold-valued smooth path as a control? Yes. It is indeed common, in physics and differential geometry, to be given a bundle $\bfB$ over some manifold $\bfM$, together with a connection, given by a $T\bfB$-valued 1-form H on $\bfM$. Lifting a smooth $\bfM$-valued path $\gamma_\bullet$ into a smooth horizontal $\bfB$-valued smooth path $\be_\bullet$ is a basic procedure in which $\gamma$ is used as a control in the ordinary differential equation 
$$
\dot \be_t = \textrm{H}(\be_t)\dot \gamma_t
$$
defining the path $\be_\bullet$. What would happen if $\gamma$ were a kind of spinning smooth $\bfM$-valued path? What about \textit{spinning geodesics} in a Riemannian setting? The notion of $p$-rough integrator introduced below, and the  results proved, will enable us to answer these questions.

\bigskip

So far there has been only a few works dealing with rough paths in a geometrical setting, starting with the work \cite{LyonsQianFlow} of Lyons and Qian. This seminal work investigated the well-posedness problem for the ordinary differential equation on the path space of a compact manifold generated by It\^o vector fields, with an eye on probabilistic applications related to path space versions of the Cameron-Martin theorem and Driver's flow equation. It was enriched by the work \cite{LiLyons} of Li and Lyons, showing that the It\^o-Lyons solution map to a Young differential equation is Fr\'echet regular under appropriate conditions, when the driving signal has finite $p$-variation, with $1<p<2$, and by the work \cite{BailleulRegularity} of the author providing a general regularity result for the It\^o-Lyons solution map, for any $p\geq 2$, in the setting of controlled paths. Another work \cite{LyonsQianJacobi} of Lyons and Qian addressed well-posedness issues for ordinary differential equations on path space associated with It\^o vector fields obtained by varying the driving rough signal.

\medskip

These works only use rough paths as an ingredient to construct some dynamics in a geometric configuration space. Cass, Litterer and Lyons made a step further in putting rough paths theory in a geometrical setting and proposed in \cite{CassLittererLyons} a notion of rough path on a manifold extending the classical notion defined in a linear setting. In the same way as a vector field on a manifold $\bfM$ can be understood in an analytic/algebraic setting as a differentiation in the ring of smooth functions on the manifold $\bfM$, a rough path is abstractly defined as a linear form on the space of sufficiently regular 1-forms on $\bfM$, which is required to have some continuity property; call it an \textit{integrator}. This functional analytic definition rests on a basic chain rule which eventually enables to understand their notion of rough path on a manifold as an equivalence class of classical rough paths, related by some chain rule under change of coordinates. This situation is exactly similar to representing a tangent vector on a $d$-dimensional manifold as an equivalence class of vectors in $\RR^d$, indexed by local diffeomorphisms of a neighbourhood of $0$ (that is local changes of coordinates), and related by a change of coordinate rule which exactly balances the changes in the numerical representation of a given 1-form $\alpha$ on $\bfM$ associated with local coordinates, so the quantity $\alpha(u)$ is independent of any choice of coordinates used to compute it. Their approach rests however on a notion of Lip-$\gamma$ manifold which prevents its easy use even with non-compact finite dimensional manifolds, not to speak about infinite dimensional manifolds.

\medskip

The ideas of \cite{CassLittererLyons} have been reloaded in a different and more accessible form in the recent work \cite{CassDriversLitterer} by Cass, Driver and Litterer, in which they define a weak geometric H\"older $p$-rough path on a finite dimensional compact embedded submanifold of $\RR^d$ as an integrator, obtained by "projection" of a weak geometric H\"older $p$-rough path in the ambiant Euclidean space, for $2\leq p<3$ only. Nothing is lost in working with submanifolds, and this notion is eventually shown to be independent of the embedding of the manifold in its environment, while these objects can be used to define and solve uniquely rough differential equations driven by weak geometric H\"older $p$-rough paths on compact manifolds. They can construct parallel transport along manifold-valued rough paths in their sense, and use it to show a one-to-one correspondence between rough paths on a finite $d$-dimensional manifold and rough paths on $d$-dimensional Euclidean space.

\medskip

We propose in this note an elementary and flexible notion of weak geometric H\"older $p$-rough path on Banach manifolds that goes beyond the previous works. Roughly speaking, a weak geometric H\"older $p$-rough path on a manifold $\bfM$ is a triple made up of a vector field valued 1-form F defined on some Banach space U, a weak geometric H\"older $p$-rough path $\bfX$ over U, and the maximal solution to the rough differential on $\bfM$ constructed from F and $\bfX$. The data of these three objects is sufficient to define and solve uniquely rough differential equations driven by "manifold-valued" rough paths, or better, by \textit{$p$-rough integrators}. The above-mentionned spinning smooth paths on $\bfM$ are precisely those $p$-rough integrators whose associated $\bfM$-valued paths are smooth. Nothing else than the (smooth) manifold structure is needed to make sense of a $p$-rough integrator. One shows in section \ref{SubsectionCanonicalRepresentation} these objects have a canonical representation if the tangent bundle of $\bfM$ is equipped with a connection, in which case one can always choose for Banach space U in the preceeding description of a $p$-rough integrator the Banach space on which the manifold is modelled.

\medskip

The interest of working with Banach manifolds comes from the fact that they naturally pop up in a number of geometric situations, as path or loop spaces over some finite dimensional manifold, as in the works \cite{BrzezniakElworthy2000, BrzezniakCarroll} of Brzezniak, Carroll and Elworthy, or the works \cite{InahamaKawabi, InahamaKawabiAbel} of Inahama and Kawabi, or as manifolds of maps of a given finite dimensional manifold, as in the works \cite{Elworthy82, BrzezniakElworthy95} of Elworthy and Brzezniak, to mention but a few works from the probability community. Note however that many interesting infinite dimensional manifolds are Fr\'echet manifolds, for which no theory of rough paths is presently available.

As a matter of fact, working with Banach manifolds will not bring any additional difficulty along the way, as compared to working with finite dimensional manifolds. We refer the reader to the books \cite{Lang} of Lang, and \cite{AbrahamMarsdenRatiu} of Abraham, Marsden and Ratiu, for the basics of differential geometry in an infinite dimensional setting. Let just recall that it makes perfect sense in a manifold setting to say that a (Banach space or real-valued) function $f$ of class $\mcC^k$, defined on the domain of some local chart $\varphi$ of $\bfM$, has bounded derivatives if the function $f\circ\varphi^{-1}$, defined on a neighbourood of $0$ in some Banach space, has bounded derivatives, as this boundedness character will not depend on the precise chart used to define it, while the bounds themselves will depend on $\varphi$.

\medskip

\noindent {\bf Notations.} We gather here a few notations that will be used throughout the note.

\begin{itemize}
   \item We shall denote by U a generic Banach space; the notation $T^{[p]}(\textrm{U})$ will be used for the truncated tensor product of order $[p]$, completed for some  choice of tensor norms. The letter $\bfX$ will stand for a weak geometric H\"older $p$-rough path over U, for some $p\geq 2$, and for the above choice of tensor norm on $T^{[p]}(\textrm{U})$.  \vspace{0.1cm} 
   \item We shall denote by $\bfM$ a Banach manifold modelled on some Banach space E. The set of continuous linear homomorphisms of E will be denoted by $\textrm{L}(\textrm{E})$, and the set of continuous linear isomorphisms of E will be denoted by $\textrm{GL}(\textrm{E})$.
\end{itemize}

\bigskip

\section{Rough integrators on Banach manifolds}
\label{SectionRPManifold}

We define in this section what will play the role in a manifold setting of weak geometric H\"older $p$-rough paths. Rough integrators are basically defined as a triple, consisting of a traditional weak geometric H\"older $p$-rough path, a vector field valued 1-form and a solution to an associated rough differential equation on the manifold. This definition, given in section \ref{SubsectionDefnRoughPaths}, requires that we recall some facts on rough differential equations with values in manifolds; this is done in section \ref{SubsectionRDEsBanach}. We explain in section \ref{SubsectionRDEs} how rough integrators can be used to define and solve rough differential equations.

\subsection{Rough differential equations}
\label{SubsectionRDEsBanach}

We adopt in this work the definition of a solution path to an $\bfM$-valued rough differential equation given in \cite{RDEsBanach} in a Banach space setting, as it is perfectly suited for our needs. It essentially amounts to requiring from a solution path that it satisfies some uniform Taylor-Euler expansion formulas, in the line of Davie' seminal work \cite{Davie}. Let U be a Banach space, and F be a $1$-form on U with values in the space of vector fields on a Banach manifold $\bfM$, of class $\mcC^{[p]+1}$. Given $u\in\textrm{U}$, we identify the vector field $\textrm{F}(\cdot\,;u)$ on $\bfM$ with its associated first order differential operator $\textrm{F}^\otimes(u)$, and extend the definition of $\textrm{F}^\otimes$ to $T^{[p]}(\textrm{U})$ setting $\textrm{F}^\otimes(1) = \textrm{Id}$, and 
$$
\textrm{F}^\otimes(u_1\otimes\cdots\otimes u_k) = \textrm{F}^\otimes(u_1)\cdots\textrm{F}^\otimes(u_k),
$$
for $1\leq k\leq [p]$, and $u_i\in \textrm{U}$, and by linearity; so $\textrm{F}^\otimes(u_1\otimes\cdots\otimes u_k)$ is a differential operator of order $k$. Given a weak geometric H\"older $p$-rough path $\bfX$ over U, recall that an $\bfM$-valued continuous path $(x_t)_{0\leq t<\zeta}$ is said to solve the rough differential equation
\begin{equation}
\label{EqManifoldValuedRDE}
dx_t = \textrm{F}^\otimes\big(x_t\,;{\bfX}(dt)\big)
\end{equation}
if there exists a constant $a>1$ such that, for any $0\leq s<\zeta$,  there exists an open neighbourhood $V_s$ of $x_s$ such that  we have the Taylor-Euler expansion
$$
f(x_t) = \big(\textrm{F}^{\otimes}\big({\bfX}_{ts}\big)f\big)(x_s) + O\big(|t-s|^a\big) 
$$
for all $t$ close enough to $s$ for $x_t$ to belong to $V_s$, and any function $f$ of class $\mcC^{[p]+1}$ defined on $V_s$, where it has bounded derivatives. The results of \cite{RDEsBanach} show that such a rough differential equation has a unique maximal solution started from any given point, as awaited. It also follows from the results of \cite{RDEsBanach}, or other classical works, that the solution path $x_\bullet$ depends continuously on the driving signal $\bfX$ in the following sense. Fix $T<\zeta$ and cover the compact support of the path $(x_t)_{0\leq t\leq T}$ by finitely many local chart domains $(O'_i)_{1\leq i\leq N}$ and $(O_i)_{1\leq i\leq N}$, with $O'_i\subset O_i$ for all $1\leq i\leq N$. Then, there exists a positive constant $\delta$ such that for $\bfY$ $\delta$-close to $\bfX$ in the H\"older rough path distance, the solution path $y_\bullet$ to the rough differential equation 
$$
dy_t = \textrm{F}^\otimes\big(y_t\,;{\bfY}(dt)\big)
$$
started from $x_0$, will be well-defined on the time interval $[0,T]$ and will remain in the open neighbourhood $\bigcup_{i=1}^N O_i$ of the support of $(x_t)_{0\leq t\leq T}$, with $y_t$ in $O_i$ whenever $x_t$ is in $O'_i$.

\ssk

Note here that the rough path setting allows to work with any Banach manifold, without any restriction on the Banach space on which it is modelled, unlike what needs to be done in a stochastic setting where a robust theory of stochastic integration is only available for some special types of Banach spaces, the so-called M-type 2 and UMD spaces, see \cite{BrzezniakElworthy,NeervenVeraarWeis}.

\bigskip

\noindent {\sc Examples.} \textit{{\bf 1. Developping a spinning straight line.} Let $\bfM$ be a 2-dimensional Riemannian manifold; denote by $H_1$ and $H_2$ the canonical horizontal vector fields on the orthonormal frame bundle $\OO\bfM$ of $\bfM$, equipped with Levi-Civita connection. Denote by $V_{12}$ the unique canonical vertical vector field and denote by $\Omega$ the (scalar-valued) curvature (tensor). Let $(\ep_1,\ep_2)$ stand for the canonical basis of $\RR^2$, and let $\bfX$ be the weak geometric H\"older $p$-rough path with first level $(t-s)\ep_1$ and area process $(t-s)\begin{pmatrix}
0 & 1 \\
-1 & 0
\end{pmatrix}$; this path was called above the \emph{spinning straight line}. Denote by \emph{H} the vector field-valued \emph{1}-form $u\in\RR^2\mapsto u^1H_1+u^2H_2$. Then the $\OO\bfM$-valued solution path to the rough differential equation
\begin{equation}
\label{EqHorizontalLift}
d\be_t = \textrm{\emph{H}}^\otimes\big(\be_t;{\bfX}(dt)\big)
\end{equation}
is actually the solution path to the ordinary differential equation
$$
\dot \be_t = \big(H_1+\Omega\,V_{12}\big)(\be_t).
$$}

\textit{{\bf 2. Finding back a geodesic.} Assume for simplicity that $\Omega$ is constant, and take as a rough path in equation \eqref{EqHorizontalLift} the weak geometric H\"older $q$-rough path $\bfY$, with $4<q<5$, whose logarithm is given by 
$$
\log {\bfY}_{ts} = (t-s)\Big(\ep_1\oplus[\ep_1,\ep_2]\oplus \Big[-\ep_1,\big[\ep_1,[\ep_1,\ep_2]\big]\Big]\Big).
$$
Then the solution path to the rough differential equation
$$
d\be_t = \textrm{\emph{H}}^\otimes\big(\be_t;{\bfY}(dt)\big)
$$
also solves the ordinary differential equation
$$
\dot \be_t = H_1(\be_t),
$$
that is, this path is the canonical lift to $\OO\bfM$ of a geodesic. This can be seen directly from the definition, or by appealing to the results proved in \cite{FrizOberhauser} by Friz and Oberhauser.}

\bigskip

\textit{{\bf 3. Rough differential equations with values in Banach Lie groups.} Let ${\sf G}$ be a Banach Lie group, with Lie algebra $\textsc{Lie}({\sf G})$. Think for instance to the loop groups, made up of $\mcC^k$ maps from some finite dimensional manifold $M_0$ to some finite dimensional Lie group $G_0$, with pointwise multiplication and inversion operations. Their Lie algebra is the set $\mcC^k\big(M_0,\frak{g}_0\big)$ of $\mcC^k$ maps from $M_0$ to the Lie algebra $\frak{g}_0$ of $G_0$, with bracket defined pointwise by the relation $\big[{\sf u},{\sf v}\big](x) := \big[{\sf u}(x),{\sf v}(x)\big]$, for any ${\sf u,v}$ in $\mcC^k\big(M_0,\frak{g}_0\big)$ and $x\in M_0$. These groups are extensively used in gauge theory or quantum field theory.}

\ssk

\textit{Write $L_g$ for the left translation by $g$ in ${\sf G}$. One defines a continuous linear map from $\textsc{Lie}({\sf G})$ to the space of smooth vector fields on ${\sf G}$ setting
$$
\textrm{\emph{F}}(g\,;u) := \Big(D_{\textrm{\emph{Id}}}L_g\Big)(u),
$$
for any $g\in {\sf G}$ and $u\in\textsc{Lie}({\sf G})$; so $\textrm{\emph{F}}(\cdot\,;u)$ is for any $u\in\textsc{Lie}({\sf G})$ a smooth left invariant vector field on ${\sf G}$. It is elementary to proceed as in the proof of theorem 4.20 in \cite{CassDriversLitterer} and see that for any weak geometric H\"older $p$-rough path $\bfX$ over $\textsc{Lie}({\sf G})$, defined on some interval $[0,T)$, for some $0\leq T\leq \infty$, the solution path to equation \eqref{EqRDE} cannot explode, so it is also defined on the interval $[0,T)$.
}

\ssk

\subsection{Rough integrators on Banach manifolds}
\label{SubsectionDefnRoughPaths}

The definition of a $p$-rough integrator given below is best understood in the light of the following example. Given a 1-form $\alpha$ on a $d$-dimensional manifold $\bfM$ equipped with a connection, how can one define the integral $\int_0^1\alpha({\circ d}x_t)$, of $\alpha$ along an $\bfM$-valued continuous semimartingale $(x_t)_{0\leq t\leq 1}$? The classical answer consists in showing that $x_\bullet$ can always be constructed as the projection in $\bfM$ of a $\textrm{GL}(\bfM)$-valued path $(\be_t)_{0\leq t\leq 1}$ obtained by solving a stochastic differential equation in $\textrm{GL}(\bfM)$ involving the horizontal vector fields associated with the connection and some $\RR^d$-valued semimartingale $(w_t)_{0\leq t\leq 1}$ -- see for instance Chapter 2 in the book \cite{Hsu} of Hsu. One then defines $\int_0^1\alpha({\circ d}x_t)$ in terms of $w$ only, setting
$$
\int_0^1\alpha({\circ d}x_t) = \int_0^1\big(\alpha\circ \be_t^{-1}\big){\circ dw_t}.
$$
So the datum of $w$ and the horizontal vector fields are all we need to define $x_\bullet$ and use it as a control. The next definition adopts a similar point of view in the present setting.

\ssk

\begin{defn}
\label{DefnRoughPathManifold}
Let $p\geq 2$ be given. A \textbf{weak geometric H\"older $p$-rough path on $\bfM$} is a a triple $\Theta = \big((x_t)_{0\leq t<\zeta},\textrm{\emph{F}},\bfX\big)$, where
\begin{itemize}
  \item $\bfX$ a weak geometric H\"older $p$-rough path over some Banach space \emph{U}, defined on the time interval $[0,\zeta)$,
  \item \emph{F} is a continuous linear map from \emph{U} to the space of vector fields on $\bfM$ of class $\mcC^{[p]+1}$,   
  \item the path $(x_t)_{0\leq t<\zeta}$ solves the rough differential equation
  \begin{equation}
  \label{EqRDE}
  dx_t = \textrm{\emph{F}}^\otimes\big(x_t\,;{\bfX}(dt)\big).  
  \end{equation}
\end{itemize}
We also call $\Theta$ a \textbf{basic $p$-rough integrator}.
\end{defn}

\ssk

Given a $p$-rough integrator $\Theta$ as above, the triple $\big((x_t)_{0\leq t\leq T},\textrm{F},\bfX\big)$, with $T<\zeta$, is also called a basic $p$-rough integrator. Let us insist on the fact that the Banach space U in the above definition is not fixed a priori, and depends on the basic $p$-rough integrator $\Theta$ on $\bfM$. So any weak geometric H\"older $p$-rough path $\bfY$ over some Banach space V can be canonically seen as a weak geometric H\"older $p$-rough path in the above sense by choosing $\textrm{U}=\textrm{V}$ and $\bfX=\bfY$, the identity map for F, and the first level of the rough path ${\bfX}_t$ for $x_t$. Note that we do not endow the set of basic $p$-rough integrators with a topology of a distance as this is somewhat subtle and we do not need it in the sequel.

\ssk

The definition of a solution to the rough differential equation \eqref{EqRDE} makes it clear that if $g$ is any sufficiently regular diffeomorphism between $\bfM$ and another manifold $\bfM'$, one defines a weak geometric H\"older $p$-rough path on $\bfM'$ setting
$$
g^*\Theta := \Big(\big(g(x_t)\big)_{0\leq t<\zeta},g^*\textrm{F},\bfX\Big),
$$
where $g^*\textrm{F}(y\,;u) := \big(D_{g^{-1}(y)}\big)g\Big(\textrm{F}\big(g^{-1}(y)\,;u\big)\Big)$, is the push forward on $\bfM'$ by $g$ of the vector field $\textrm{F}(\cdot\,;u)$ on $\bfM$.

\medskip

Definition \ref{DefnRoughPathManifold} encompasses the characterization of a weak geometric H\"older $p$-rough path on a submanifold of $\RR^d$ given in \cite{CassDriversLitterer}, only for $2\leq p<3$, in terms of a priori extrinsic considerations. Roughly speaking, they define their class of weak geometric $p$-rough paths on a compact submanifold $\bfM$ of $\RR^d$, for $2\leq p<3$, as the class of all weak geometric rough paths $\bfY$ over $\RR^d$ with the property that ${\bfY}_0=(1,y_0,0,\dots)$ with $y_0\in\bfM$, and
\begin{equation}
\label{EqRPCassDriverLitterer}
d{\bfY}_t = Q(y_t)\,d{\bfY}_t
\end{equation}
 for some projector-valued sufficiently regular map $Q$ such that $Q(y)$ has range included in $T_y{\bfM}$ if $y\in\bfM$; this equation ensures that $y_t\in\bfM$. If ${\bfM}'\subset\RR^{d'} $ is any other embedded compact manifold diffeomorphic to $\bfM$ through $\varphi : \bfM \rightarrow \bfM'$, equation \eqref{EqRPCassDriverLitterer} and the chain rule for rough integrals shows that one defines an image rough path $\bfZ$ over $\RR^{d'}$ setting
 $$
 d{\bfZ}_t = d\varphi(d{\bfY}_t) 
 $$ 
and that it satisfies the identity
$$
\int_0^1 (\varphi^*\alpha)(d{\bfZ}_t) = \int_0^1 \alpha(d{\bfY}_t),
$$
for any sufficiently regular 1-form $\alpha$ on $T\bfM$. (The expression $\varphi^*\alpha$ stands for the push-forward of the 1-form $\alpha$ by $\varphi$, and both integrals make sense because of condition \eqref{EqRPCassDriverLitterer}.) So the class of weak geometric $p$-rough paths is intrinsically defined as a class of rough integrators, independently of any particular embedding of $\bfM$, seen as an abstract manifold. Note that they use general controls $\omega(s,t)$ to measure the size of their rough paths, while we simply use $\omega(s,t)=t-s$ and H\"older scales. A weak geometric $p$-rough path in the sense of \cite{CassDriversLitterer} is simply described in our terms as a special kind of weak geometric $p$-rough path in the ambiant linear space via equation \eqref{EqRPCassDriverLitterer}. The interest of working with the intrinsic definition \ref{DefnRoughPathManifold} is clear in an infinite dimensional setting where embeddings of Banach manifolds into larger Banach spaces are less natural and rarely happen, unless $\bfM$ is parallelizable \cite{EellsElworthy}. Note however that the frame bundle of $\bfM$ is always parallelizable if $T\bfM$ is equipped with a connection. See section \ref{SectionCanonicalRepresentation}.

\medskip

One of the nice features of the notion of rough integrator on a manifold proposed in \cite{CassLittererLyons} and \cite{CassDriversLitterer} is the possibility to push forward a rough integrator by a sufficiently regular map $g$ from $\bfM$ to another manifold $\bfN$, giving as a result a rough integrator on $\bfN$. The rough integral of a 1-form $\alpha$ on $\bfN$ against the image rough integrator is simply defined as the rough integral of $\alpha\circ dg$ against the rough integrator on $\bfM$. A similar transport operation can be done in our framework, by seeing the push forward of a $p$-rough integrator $\Theta$ on $\bfM$ by $g : \bfM\rightarrow\bfN$ as the $p$-rough integrator $g\Theta$ on $\bfM\times\bfN$ associated with $(\textrm{F},dg\circ\textrm{F})$ and $\bfX$.

\medskip

It may happen that the vector fields $\textrm{F}(\cdot\,;u)$ and $\textrm{F}(\cdot\,;v)$ commute for any $u,v\in \textrm{U}$, and $\bfX$ has null first level, so the path $x_\bullet$ is constant. We say in that case that $\Theta$ is a \textbf{pure rough path}. It can still give rise to some dynamics, when seen as a rough integrator in a rough differential equation, as shown in the next section.

\subsection{Weak geometric H\"older $p$-rough paths as integrators}
\label{SubsectionRDEs}

Recall the construction of a line integral along a continuous semimartingale on a manifold described at the beginning of section \ref{SubsectionDefnRoughPaths}. The next proposition shows that a weak geometric H\"older $p$-rough path on $\bfM$ can be used as an integrator in a rough differential equation, which justifies calling it a \textit{$p$-rough integrator}. 

\begin{defn}
\label{DefnRDESolutionManifold}
Let $\Theta$ be a weak geometric H\"older $p$-rough path on $\bfM$ in the sense of definition \emph{\ref{DefnRoughPathManifold}}. Let also $\pi : \bfB\mapsto \bfM$ be a smooth fiber bundle over $\bfM$, equipped with a connection, of class $\mcC^{[p]+1}$, given under the form of a smooth $T\bfB$-valued \emph{1}-form on $\bfM$ such that $\pi_* : T_\be\bfB \rightarrow T_{\pi(\be)}\bfM$ is an isomorphism for every $\be\in\bfB$. Given a rough integrator $\Theta = \big((x_t)_{0\leq t<\zeta},\textrm{\emph{F}},\bfX\big)$ in $\bfM$, and $\be_0\in \bfB$ with $\pi(\be_0)=x_0$, a path $(\be_t)_{0\leq t< \zeta'}$ in $\bfB$ is said to be a lift of $(x_t)_{0\leq t<\zeta'}$, if $\zeta'\leq \zeta$ and it is the solution path to the rough differential equation
\begin{equation}
\label{EqRDEManifoldManifold}
d\be_t = \textrm{\emph{H}}\Big(\textrm{\emph{F}}\big(\pi(\be);{\bfX}(dt)\big)\Big)
\end{equation}
started from $\be_0$. This equation ensures in particular that $\pi(\be_t)=x_t$. The $p$-rough integrator $\big((\be_t)_{0\leq t<\zeta'},\textrm{\emph{H}}\circ\textrm{\emph{F}}\circ\pi,\bfX\big)$ is then said to be the \textbf{\emph{lift of $\Theta$ to $\bfB$.}} 
\end{defn}

A commonly encountered situation is to have $\bfB = \bfM\times\bfN$, for some other Banach manifold $\bfN$, and a connection form defined for any $p\in T_x\bfM$ and $y\in\bfN$ by the formula
$$
\textrm{H}(x,y)p = \big(p,\textrm{G}(y,x;p)\big),
$$
with $T(\bfM\times\bfN)$ canonically identified with $T{\bfM}\times T\bfN$, and some function $\textrm{G} : {\bfN}\times T{\bfM}\rightarrow T\bfN$ of class $\mcC^{[P]+1}$, with $\textrm{G}(y,x;p)\in T_y\bfN$, for each $p\in T_x\bfM$, and $x\in\bfM$.

\ssk

This definition is reminiscent of the approach chosen by Cass, Litterer and Lyons, in section 5 of \cite{CassLittererLyons}, to define in their setting rough differential equations driven by manifold-valued rough paths. Note however that their theory requires the use of the rigid and non-trivial notion of $\textrm{Lip}-\gamma$ manifold; no such constraint holds in our setting.

\medskip

The fact that the two vector fields $\big(\textrm{H}\circ\textrm{F}\circ\pi\big)(\cdot,\cdot\,;u)$ and $\big(\textrm{H}\circ\textrm{F}\circ\pi\big)(\cdot,\cdot\,;v)$ may not commute while $\textrm{F}(\cdot\,;u)$ and $\textrm{F}(\cdot\,;v)$ may commute, for some $u,v\in\textrm{U}$, explains why pure rough paths can generate dynamics. The results on rough differential equations recalled in the introduction of this section apply and show that 

\begin{prop}
The rough differential equation \eqref{EqRDEManifoldManifold} has a unique maximal solution started from any given point.
\end{prop}

Solving successively some rough differential equations of the form \eqref{EqRDEManifoldManifold} with some basic $p$-rough integrators $\Theta^{(1)}=\big(\big(x_t^{(1)}\big)_{0\leq t\leq t_1},\textrm{F}^{(1)},{\bfX}^{(1)}\big),\dots,\Theta^{(k)}$, with $x_0^{(j)}=x_{t_{j-1}}^{(j-1)}, \;y_0^{(j)}=y_{t_{j-1}}^{(j-1)}$ and $t_{j-1}<\zeta'_{j-1}\leq \zeta_{j-1}$, for $2\leq j\leq k$, defines, whenever this makes sense, the \textit{concatenation $\Theta^{(k)}\cdots\Theta^{(1)}$ of the basic $p$-rough integrators $\Theta^{(i)}$}. We call such an object a {\bf $p$-rough integrator}.

\medskip

\section{Canonical representation of rough integrators}
\label{SectionCanonicalRepresentation}

We show in this section that $p$-rough integrators have a canonical representation when the tangent bundle of $\bfM$ is equipped with a connection. This representation is the analogue of the representation of a regular path $\gamma$ on $\bfM$ by a regular path in $T_{\gamma_0}\bfM$ using Cartan's development map.

\subsection{Cartan's moving frame method}
\label{SubsectionGeometry}

Let $\bfM$ be a manifold of finite dimension $d$. One owes to E. Cartan the introduction in differential geometry of the \textit{moving frame method}, which provides a chart for the set of $\bfM$-valued paths of class $\mcC^1$ based at some fixed starting point, in terms of $\mcC^1$ paths in $\RR^d$ started from $0$. Its construction requires the use of a connection on $T\bfM$ and of the frame bundle of $\bfM$; it can basically be described as follows. Given a $\mcC^1$ path $\gamma_\bullet = (\gamma_t)_{0\leq t\leq 1}$ on $\bfM$ and a frame $\be_0$ at $\gamma_0$, its parallel transport along the path $\gamma_\bullet$ defines a section $(\be_t)_{0\leq t\leq 1}$ of the frame bundle above $\gamma_\bullet$. One can use these frames to describe $\dot \gamma_t$ at any time in terms of its coordinates $\big(\dot u^i_t\big)_{1\leq i\leq d}$ in $\be_t$. The $\RR^d$-valued path $u_\bullet$ defined  by the formula
$$
u_t: = \int_0^t \dot u^i_s\,ds
$$
is called the \textbf{anti-development of the path} $\gamma$. One finds back $\gamma_\bullet$ from $u_\bullet$ and $(x_0,\be_0)$ as the projection on $\bfM$ of the unique solution of the ordinary differential equation in the frame bundle
\begin{equation*}
\nabla_{\dot \gamma_t}\be_t = 0, \quad \dot \gamma_t := \be_t\big(\dot u_t\big).
\end{equation*}
The interest of Cartan's moving frame method is that it somehow provides a dimensionally-optimal coding of a path in $\bfM$ in terms of vector space-valued paths. Here is a trivial illustration of this fact. Let $u_\bullet$ be an $\RR^{10}$-valued path and $\gamma_t = \sum_{i=1}^{10} u^i_t$ be an $\RR$-valued path. The coding of $\gamma_\bullet$ by the 10 components of $u_\bullet$ is optimized by coding it with its real value at each time. Replacing $\RR^{10}$ and $\RR$ by an infinite dimensional spaces emphasizes the importance of such parcimonious representations.

\bigskip

Recall $\bfM$ is modelled on some Banach space E. One can actually give a parcimonious description of weak geometric H\"older $p$-rough paths on Banach manifolds similar to the above one, providing a description of these objects in terms of weak geometric H\"older $p$-rough paths in the tensor space $T^{[p]}(\textrm{E})$, as opposed to the a priori unrelated Banach space U involved in the primary definition of a weak geometric H\"older $p$-rough path, see definition \ref{DefnRoughPathManifold}. As in Cartan's moving frame method, this requires the tangent bundle $T\bfM\rightarrow\bfM$ to be equipped with a connection; this is a non-trivial assumption in an infinite dimensional setting, linked to the fact that there exists some Banach spaces that do not even admit a smooth partition of unity. No such pathology happens in finite dimension, and finite dimensional manifolds can always be endowed with an arbitrary connection; so the results of the forthcoming section \ref{SubsectionCanonicalRepresentation} always hold for rough integrators on finite dimensional manifolds.

\medskip

The frame bundle $\textrm{GL}(\bfM)$  of $\bfM$ will play a crucial role in that play. This is the collection of all isomorphisms from E to $T_m\bfM$, for $m\in \bfM$; it has a natural manifold structure modelled on $\textrm{GL}(\textrm{E})\times\textrm{E}$. We shall denote by $\be$ a generic element of $\textrm{GL}(\bfM)$, and by $\mcV\textrm{GL}(\bfM)$ the vertical sub-bundle of $T\textrm{GL}(\bfM)$, canonically identified with the Lie algebra $\textrm{gl}(\textrm{E}) = \textrm{L}(\textrm{E})$ of $\textrm{GL}(\textrm{E})$. The connection $\nabla$ on the bundle $T\bfM\rightarrow \bfM$ is naturally lifted into a connection on the bundle $\pi : \textrm{GL}(\bfM)\rightarrow \bfM$, still denoted by the same symbol. Remark that $\pi_*$ is an isomorphism between the horizontal sub-bundle $\mcH^\nabla\textrm{GL}(\bfM)$ in $T\textrm{GL}(\bfM)$ and $T\bfM$. The horizontal destribution in $T\textrm{GL}(\bfM)$ can be used to define a continuous linear map ${\bfH}^\nabla$ from E to the space of horizontal vector fields on $T\textrm{GL}(\bfM)$, defined by the requirement that ${\bfH}^\nabla(\be\,;a)\in T_\be\textrm{GL}(\bfM)$ is horizontal and corresponds to $\be(a)\in T\bfM$, for any $a\in \textrm{E}$ and $\be\in \textrm{GL}(\bfM)$. This vector field valued 1-form ${\bfH}^\nabla$ on E is called the \textbf{canonical horizontal 1-form}.

\ssk

Once again, we refer the reader to the nice books \cite{Lang} and \cite{AbrahamMarsdenRatiu} for the basics of differential geometry on Banach manifolds.

\ssk

\subsection{Canonical representation of rough integrators}
\label{SubsectionCanonicalRepresentation}

Theorem \ref{ThmCanonicalRepresentation} below gives a canonical and parcimonious representation of a given weak geometric H\"older $p$-rough path, seen as a rough integrator. We assume for that purpose that the manifold $\bfM$ is endowed with a connection $\nabla$. Given a 1-form F on U with values in the space of vector fields on $\bfM$ as in the above definition of a weak geometric H\"older $p$-rough path on $\bfM$, we denote by $\bfF$ its lift to a 1-form with values in the space of $\nabla$-horizontal vector fields on $\textrm{GL}(\bfM)$. Parallel translation along a weak geometric H\"older $p$-rough path $\Theta = \big((x_t)_{0\leq t<\zeta},\textrm{F},\bfX\big)$ is defined as the solution path to the rough differential equation in $\textrm{GL}(\bfM)$
$$
d\be_t = {\bfF}\big(\be_t,{\bfX}(dt)\big)
$$
started from any given frame $\be_0\in\textrm{GL}(\bfM)$ above $x_0$. It is elementary to see that the paths $\be_\bullet$ and $x_\bullet$ are defined on the same maximal interval.

\ssk

Given another Banach manifold $\bfN$, and a 1-form $\textrm{G} : {\bfN}\times T{\bfM}\rightarrow T\bfN$, as used above to write down a rough differential equation in definition \ref{DefnRDESolutionManifold}, define $\frak{S}$ on $\textrm{GL}({\bfM})\times {\bfN}\times \textrm{U}$ by the formula
$$
\frak{S}(\be,y\,; u) = \Big({\bfF}(\be\,;u),\,\textrm{G}\big(y,\pi(\be)\,;\pi_*{\bfF}(\be\,;u)\big)\Big),
$$
and $\overline{\frak{S}}$ on $\textrm{GL}({\bfM})\times {\bfN}\times \textrm{E}$ by the formula
$$
\overline{\frak{S}}(\be,y\,; a) = \Big({\bfH}^\nabla(\be\,;a),\,\textrm{G}\big(y,\pi(\be)\,;\pi_*{\bfH}^\nabla(\be\,;a)\big)\Big).
$$

\ssk

\begin{thm}[Canonical representation of rough integrators]
\label{ThmCanonicalRepresentation}
Let $\bfM$ be a Banach manifold endowed with a connection $\nabla$, and $\Theta = \big((x_t)_{0\leq t<\zeta},\textrm{\emph{F}},\bfX\big)$ be a basic $p$-rough integrator on $\bfM$.
\begin{enumerate}
   \item One defines a weak geometric H\"older $p$-rough path $\bfZ$ over \emph{E}, on the time interval $[0,\zeta)$, by solving the rough differential equation
\begin{equation}
\label{EqConstructionZ}
\begin{split}
&d\be_t = {\bfF}\big(\be_t\,;{\bfX}(dt)\big),  \\
&d{\bfZ}_t = {\bfZ}_t\otimes \be_t^{-1}\textrm{\emph{F}}\big(\pi(\be_t)\,;{\bfX}(dt)\big). \vspace{0.1cm}
\end{split}
\end{equation}
in $\textrm{\emph{GL}}({\bfM})\times T^{[p]}(\textrm{\emph{E}})$, started from $(\be_0,\textrm{\emph{Id}})$.  \vspace{0.1cm}
   \item The solution paths $(\be_t,y_t)$ and $(\overline{\be}_t,\overline{y}_t)$ in $\textrm{\emph{GL}}({\bfM})\times\bfN$ to the rough differential equations 
   \begin{equation}
   \begin{split}
   &d(\be_t,y_t) = \frak{S}\big(\be_t,y_t\,;{\bfX}(dt)\big), \\
   &d(\overline{\be}_t,\overline{y}_t) = \overline{\frak{S}}\big(\overline{\be}_t,\overline{y}_t\,;{\bfZ}(dt)\big),   
   \end{split}
   \end{equation}
   coincide if they start from the same initial condition.
\end{enumerate}
\end{thm}

\ssk

\noindent It is elementary to state an analogue of this theorem when $\bfM\times\bfN$ is replaced by a fiber bundle $\bfB$ over $\bfM$, equipped with a regular connection. There is no  more difficulties in stating a similar result for any $p$-rough integrator instead of basic $p$-rough integrators. The interesting point is that although the rough paths in each $\Theta^{(i)}$ may be defined on {\it different}  Banach spaces, the canonical representation of the concatenation of the $\Theta^{(i)}$ only involves a rough path over E.

\medskip

\begin{Dem}
The first claim follows from general principles. For the second point, pick $p<p'<[p]+1$. Using the continuity result for the It\^o map recalled in section \ref{SubsectionRDEsBanach}, in the topology of H\"older $p'$-rough path, together with the continuous embedding of the space $WG_p(\textrm{U})$ of weak geometric H\"older $p$-rough paths in U into the space $G_{p'}(\textrm{U})$ of geometric H\"older $p'$-rough paths in U, and the continuous embedding of $G_{p'}(\textrm{U})$ into $WG_{p'}(\textrm{U})$, we see that it suffices to prove the claim when $\bfX$ is the weak geometric H\"older ($p'$-) rough path associated with a smooth U-valued control. The result is clear in that case as it appears as a rephrasing of Cartan's moving frame method, as described above in section \ref{SubsectionGeometry}.
\end{Dem}

\bigskip

So one can always understand the solution of a rough differential equation driven by a $p$-rough integrator $\Theta$ as the solution to another rough differential equation driven by a $p$-rough integrator involving a weak geometric H\"older $p$-rough path over the model space E, and the canonical horizontal 1-form $\bfH^\nabla$. The dependence on F in the original $p$-rough integrator is hidden in the definition fo $\bfZ$ in this reformulation. Given $\be_0$ above $x_0$, the "$\textrm{GL}(\bfM)$-valued" $p$-rough integrator $\big((\be_t)_{0\leq t<\zeta},{\bfH}^\nabla,\bfZ\big)$ is said to be the \textbf{canonical representation of $\Theta$}.

\medskip

This result echoes one of the main results of \cite{CassDriversLitterer}, corollary 6.19, where a one-to-one correspondence between the classes of weak geometric rough path on a given finite $d$-dimensional manifold, weak geometric horizontal rough paths on its frame bundle, and classical weak geometric rough paths on $\RR^d$, is proved in their setting, using Cartan's development map as an essential ingredient. Applications of our setting will be developped in a forthcoming work.

\end{document}